\documentclass[11pt]{article}
\usepackage{graphicx}
\usepackage{amssymb}
\usepackage{amsmath}

\newcommand{\ra}{\rangle}
\newcommand{\la}{\langle}
\newcommand{\goesto}{\longrightarrow}
\newcommand{\dif}[2]{\frac{\partial #1}{\partial#2}}
\newcommand{\dv}[1]{\text{ div}{\left(#1\right)}}
\newcommand{\R}{\mathbb{R}}

\newcommand{\setin}{\subset}

\newtheorem{lemma}{Lemma}
\newtheorem{theorem}{Theorem}
\newtheorem{corollary}[theorem]{Corollary}

\title{An Inviscid Regularization for
the Surface Quasi-Geostrophic Equation }
\author{Boualem Khouider\thanks{Department of Mathematics and Statistics, University of Victoria, PO BOX 3045 STN CSC Victoria, B.C., Canada, V8W 3P4, Email: khouider@math.uvic.ca} ~ and Edriss S. Titi\thanks{Department of Mathematics and  Department of Mechanical and  Aerospace Engineering, University of California Irvine, CA  92697-3875, USA {\bf also} Department of Computer Science and Applied Mathematics, Weizmann Institute of Science Rehovot 76100, Israel}}
\date{January 30, 2007}
\begin{document}
\maketitle

{\bf Abstract}\\
 Inspired by recent developments  in Berdina-like models
 for turbulence,  we propose an inviscid regularization
 for the surface quasi-geostrophic (SQG) equations.
 We are particularly interested in  the celebrated
 question of blowup in finite time of the solution
 gradient of the SQG equations. The new regularization
 yields a necessary and sufficient condition, satisfied
 by the regularized solution, when a regularization
 parameter $\alpha$ tends to zero, for the solution of
 the original SQG equations to develop a singularity in
 finite time.  As opposed to the commonly used viscous
 regularization,  the inviscid  equations derived here
 conserve a modified energy. Therefore, the new regularization
 provides an attractive numerical procedure for finite time
 blow up testing. In particular, we prove that, if the initial
 condition is smooth, then the regularized solution remains as
 smooth as the initial data for all times. Moreover, much like
 the original problem, the inviscid regularization admits a
 maximum principle.

\section{Introduction}
We consider the surface quasi-geostrophic (SQG) equations,  with periodic boundary conditions on a basic periodic square $\Omega = [0,1]^2 \setin{\R}^2$,
\begin{align}
\dif{\theta}{t} + \text{div}({\bf v}\theta) = 0 \notag\\
(-\Delta)^{1/2} \psi = \theta \label{sqg}\\
 \nabla^\bot  \psi = {\bf v}\notag, \notag\\
 \theta(x,0) = \theta_0(x) & \text{ in } \Omega \notag\\
\int_\Omega \theta \,d{\bf x}=0,\int_\Omega \psi \,d{\bf x}=0,\int_{\Omega} {\bf v} \,d{\bf x}=0 \notag .
\end{align}
Here $\displaystyle \Delta = \frac{\partial^2}{\partial x^2} + \frac{\partial^2}{\partial y^2}$ is the horizontal Laplacian operator and $(-\Delta)^{1/2}$ is the pseudo-differential operator defined  in the Fourier space by $\widehat{ (-\Delta)^{1/2} u}({\bf k} ) = |{\bf k}| \hat u({\bf k}) $. The first equation in (\ref{sqg}) describes the evolution of the potential temperature, $\theta$, at the surface  of the   ocean, for a quasi-geostrophic flow, i.e. a first order perturbation of a ''geosptrophically balanced'' mean state; a state where  the horizontal pressure gradient is balanced by the vertical  component of the Coriolis force. Details on the derivation of,  and more discussion on, the system in (\ref{sqg}) can be found in \cite{CMT,Pedloski} and references therein. In (\ref{sqg}), ${\bf v}$ represents the incompressible horizontal velocity at the surface and $\psi$ is the stream function.

The system of equations in (\ref{sqg}) is  interesting in itself since it models an important geophysical problem.  However, it has also been the focus of interesting mathematical work \cite{canswu,CF,CMT,MT}, since the evolution of $\nabla \theta$   resembles the evolution of the vorticity in the 3D Euler equations.  This is despite the 2D nature of the equation, and the misleading impression that $\theta$ evolves like the vorticity in the 2D Euler equations. Preliminary numerical simulations conducted in  \cite{MT} revealed that the SQG equations   (\ref{sqg}) with smooth initial data develop  sharp fronts in the level contours of $\theta$ and conjectured the possibility of formation of a finite time singularity in $\nabla\theta$.  A more careful simulation that was conducted in \cite{CNT1,CNT2} revealed the absence of such singularity and attributed the observation in the simualtion of \cite{MT} to a growth of the type $e^{e^t}$ .  Indeed it was proven rigorously in \cite{cord} that the scenario of blow up suggested  in \cite{MT} is not possible.

We present here a new inviscid regularization (\ref{pb}), inspired by the inviscid simplified Bardina model of turbulence \cite{cao}  (see also the inviscid version of the Navier-Stokes-Voight model \cite{oskolkov}), for the SQG equations (\ref{sqg}). This new regularizations
 yields a necessary and sufficient condition, satisfied by the regularized solution, for  $\nabla \theta$ to blow up as the regularization parameter tends to zero.  As opposed to   viscous regularizations, used extensively in  analytical studies  (see, e.g.,   \cite{canswu}), the inviscid regularization employed here conserves a modified energy (\ref{eqe}). In fact, instead of smoothing the solution by dissipating energy at small scales, here the small scales are simply prevented from getting too much energy via a penalty method due to the modified energy (see (\ref{eqe}) below). Therefore, the new regularization  provides a systematic practical procedure for finite time blow up testing for (\ref{sqg}). The rest of the paper is organized as follows. The regularized problem is introduced in section 2 where minimal regularity requirements guaranteeing   the conservation of the modified energy (\ref{eqe}).  In section 3 we prove long time existence and uniqueness for the regularized problem (\ref{pb}). In section 4 we prove that if the initial data is smooth (in $H^m$) then the solution for the regularized problem (\ref{pb}) remains smooth   (in $H^m$) for all time. Moreover, we prove that, as the original SQG-equation (\ref{sqg}), the regularized problem admits a maximum principle. In section 5, we prove that the regularized solution of (\ref{pb}) converges  to a weak solution of the SQG equations (\ref{sqg}), when the regularizing parameter $\alpha$ goes to zero.  We also prove that if the original SQG equations have a regular (smooth) solution then the regularized solution necessarily converges strongly to this regular solution. We finally, prove in section 5 that a necessary and sufficient condition, satisfied by the regularized solution $\theta^{\alpha}$ of (\ref{pb}), for the solution of the SQG equations to blow up in finite time $T$ is \[ \sup_{t\in [0,T)}\liminf_{\alpha\goesto 0}
\alpha ||\nabla \theta^{\alpha}(t)||_{L^2} = \epsilon >0,\]
where $[0,T)$ is the maximal interval of existence of solutions of (\ref{sqg}). 

Numerical tests of this approach will be reported in a forthcoming work.



\section{The regularized problem}
Let $\alpha$ be a small positive parameter (length scale). Consider the following inviscid regularization of (\ref{sqg}) \begin{align}
&(1-\alpha^2 \Delta) \theta^{\alpha} = \tilde \theta^{\alpha} \notag\\
&\dif{\tilde \theta^{\alpha}}{t} + \dv {{\bf v}^{\alpha} \theta^{\alpha}}= 0 & \text{ in } \Omega \times [0,T]\label{pb}\\
&(-\Delta)^{1/2} \psi^{\alpha} = \notag \theta^{\alpha} \\
&\nabla^\bot \psi^{\alpha} = {\bf v}^{\alpha}; \int_{\Omega}\theta^{\alpha}d\,{\bf x}=0 ; \int_{\Omega}\psi^{\alpha}\,d{\bf x}=0 \notag\\
& \theta^{\alpha}(x,0) = \theta_0(x) & \text{ in } \Omega \notag
\end{align}
subject to periodic boundary conditions. Unless otherwise stated, the super-script $\alpha$ is dropped below in sections 2-4 to simplify our presentation.

\subsection{Energy conservation and minimal regularity}
We define the modified energy for the solution $\theta$ of the regularized problem in (\ref{pb}) as \begin{equation}\label{eqe} E(t) = \int_\Omega \left(\theta^2({\bf x},t) + \alpha^2 |\nabla \theta({\bf x},t)|^2 \right)\,d{\bf x}. \end{equation}
It is easy to show that if the solution $\theta,{\bf v},\psi$ for the regularized problem (\ref{pb}) is smooth then  the energy $E(t)$ of the system is conserved. This statement will be made rigorous below in the proof of Theorem \ref{theorem1}.
Moreover, we notice that if $\theta$ is in the Sobolev space $H^1(\Omega) = \{ u\in L^{2}(\Omega): \int_{\Omega} \left(u^2 + |\nabla u|^2 \right)\,d{\bf x} <+\infty\}$,
then $\psi = (-\Delta)^{-1/2}\theta$ belongs to the Sobolev space $H^2= \{ u \in H^1(\Omega): \int_{\Omega}  |\partial_{x_i,x_j} u|^2 \,d{\bf x} <+\infty\}$,
where $\partial_{x_i,x_j}$ is any derivative of second order. This in turn implies that ${\bf v}=\nabla^\bot \psi$ is in $H^1(\Omega)$.
Therefore, under the periodic  boundary conditions, the integral
\begin{equation}\label{E1} \int_{\Omega} \text{div}({\bf v}\theta) \theta \,d{\bf x} = \frac{1}{2}\int\text{div}({\bf v}\theta^2)\,d{\bf x}=0\end{equation}
 which implies that \[0=\frac{d}{dt}\int
\tilde\theta({\bf x},t) \theta({\bf x},t)  \,d{\bf x} = \frac{d}{dt}\int\left(\theta^2 ({\bf x},t) + \alpha^2|\nabla \theta({\bf x},t) |^2\right)  \,d{\bf x},\]
where $\tilde \theta = \theta -\alpha^2\Delta \theta$.
This is made rigorous in the following lemma.

\begin{lemma} \label{lemma1}
Let  $\theta\in\dot H^1(\Omega)=\{ u \in H^1(\Omega), \int_{\Omega}u \,d{\bf x} =0\}$  and ${\bf v}\in \dot H^1 (\Omega)\times \dot H^1(\Omega)$
then $\dv{{\bf v}\theta}\in \dot H^{-1}(\Omega)$, where $\dot H^{-1}(\Omega)$ is the dual space of $\dot H^1(\Omega)$. Moreover, for ${\bf v}\in  \dot H^1 (\Omega)\times \dot H^1(\Omega)$ fixed, $\theta \longrightarrow \dv{\theta{\bf v}}$ is a linear continuous operator from $ \dot H^1 (\Omega)$ to $ \dot H^{-1}(\Omega)$.
\end{lemma}

{\em Proof:}
Let $\phi \in \dot H^1(\Omega)$, then
\begin{align} | \langle \text{div}({\bf v}\theta),\phi \rangle| &= |\int_{\Omega} \theta({\bf x}){\bf v}({\bf x})\cdot\nabla\phi ({\bf x})\,d{\bf x} | \notag \\
&\le C ||{\bf v}||_{L^4(\Omega)} ||\theta||_{L^4(\Omega)} ||{\nabla \phi}||_{L^2(\Omega)}\notag \text{ (by H\"older's inequality)}\\
 &\le C ||{\bf v}||^{1/2}_{L^2(\Omega)}||\nabla{\bf v}||^{1/2}_{L^2(\Omega)} ||\theta||^{1/2}_{L^2(\Omega)}||\nabla\theta||^{1/2}_{L^2(\Omega)} ||{\nabla \phi}||_{L^2(\Omega)},\notag\end{align} by the 2D Gagliardo-Nirenberg-Ladyzhenskaya interpolation inequality (see, e.g., \cite{Adams,CF,L,LZ2}).
Therefore\begin{align}
 &||\text{div}({\bf v}\theta) ||_{\dot  H^{-1}} \le C ||{\bf v}||^{1/2}_{L^2(\Omega)}||\nabla{\bf v}||^{1/2}_{L^2(\Omega)} ||\theta||^{1/2}_{L^2(\Omega)}||\nabla\theta||^{1/2}_{L^2(\Omega)}.\notag
\end{align}
As a consequence we have the following corollary.

\begin{corollary}
Let ${\bf v}\in \dot H^1(\Omega)\times \dot H^1(\Omega)$ such that $\text{div }{\bf v}=0$ and $\theta\in\dot H^1({\Omega})$. Then\begin{equation}\label{eq5}
\la \dv{{\bf v}\theta},\theta\ra =0.\end{equation}
\label{co1}
\end{corollary}
\section{Global existence for the regularized problem}
Here we prove that the regularized problem in (\ref{pb}) admits a global smooth solution for all time if the initial condition $\theta_0$ is smooth. More precisely, we have the following theorem.
\begin{theorem}\label{theorem1}
Let $\theta_0 \in \dot H^1(\Omega)$ (i.e. $\tilde \theta_0 = (1-\alpha^2\Delta)\theta_0 \in \dot  H^{-1}(\Omega)$), then the initial value problem
\begin{align}&(1-\alpha^2 \Delta) \theta= \tilde \theta \notag\\
& \partial_t \tilde \theta =- \text{\em div}(\theta {\bf v}),\notag\\
& {\bf v}= \nabla^\bot \psi, (-\Delta)^{1/2} \psi = \theta, \int_{\Omega}\psi \,d{\bf x} =0 \label{eq3}\\
&\theta(0)= \theta_0 \notag \end{align}
has a global unique solution $\theta \in \mathcal{C}^1((-\infty,+\infty),\dot H^1(\Omega))$ (or, $\tilde \theta \in \mathcal{C}^1((-\infty,+\infty),\dot H^{-1}(\Omega))$).
\end{theorem}
{\em Proof:} \\
Given the relations $\tilde \theta = (1-\alpha^2\Delta)\theta, {\bf v}= \nabla^\bot \psi, (-\Delta)^{1/2} \psi = \theta$, we can write
\[- \text{div}(\theta{\bf v}) = F(\tilde \theta):= -\text{div}\left( \left(\nabla^\bot ( -\Delta)^{-1/2}(1-\alpha^2\Delta)^{-1}\tilde \theta\right)  \times (1-\alpha^2\Delta)^{-1}\tilde\theta\right).\]
Thus, by Lemma \ref{lemma1},  we have a functional differential equation of the form
\begin{equation}\label{fde} \frac{d}{dt}\tilde\theta = F(\tilde \theta), ~\tilde \theta(0)=\tilde \theta_0\end{equation}
in the space $\dot H^{-1}(\Omega)$.
We first show short time existence and uniqueness. For this, it is enough to establish that the functional $F(\tilde \theta)$ is locally Lipshitz as a map from $\dot H^{-1}(\Omega)$ into $\dot H^{-1}(\Omega)$.
We have
\begin{align} ||F(\tilde \theta_1)-F(\tilde \theta_2)||_{\dot H^{-1}} &= || \text{div}(\theta_1 {\bf v}_1)-\text{div}(\theta_2 {\bf v}_2)||_{\dot H^{-1}}
\notag\\
&\le ||\text{div}(\theta_1 ({\bf v}_1-{\bf v}_2))||_{\dot H^{-1}} +
||\text{div}((\theta_1 -\theta_2){\bf v}_2)||_{\dot H^{-1}}\notag \\
&\le C\Big( ||{\bf v}_1-{\bf v}_2||^{1/2}_{L^2}||\nabla({\bf v}_1-{\bf v}_2)||^{1/2}_{L^2} ||\theta_1||^{1/2}_{L^2}||\nabla\theta_1||^{1/2}_{L^2} \notag \\ &~~~~~~~~~~+ ||{\bf v}_2||^{1/2}_{L^2}||\nabla{\bf v}_2||^{1/2}_{L^2} ||\theta_1-\theta_2||^{1/2}_{L^2}||\nabla(\theta_1-\theta_2)||^{1/2}_{L^2}\Big)\notag
\end{align}
Now we invoke Poincar\'e inequality \cite{Adams,CF,LZ2,MB}
\begin{equation}\label{PI} ||\phi||_{L^2} \le \lambda_1^{-1/2} ||\nabla \phi||_{L^2}, \forall \phi \in \dot H^{1}(\Omega).\end{equation}
Here $\lambda_1$ is the first eigenvalue of $-\Delta$ with domain $D(-\Delta)= H^2(\Omega)\cap \dot H^(\Omega)$.
Which leads to
\[ ||F(\tilde \theta_1) -F(\tilde \theta_2)||_{\dot H^{-1} }\le C ( ||\nabla({\bf v}_1-{\bf v}_2)||_{L^2} ||\nabla\theta_1||^{1/2}_{L^2} +  ||\nabla{\bf v}_2||_{L^2} ||\nabla(\theta_1-\theta_2)||_{L^2}
).\]

But, given that the functional operator $\theta\longrightarrow {\bf v} = \nabla^\bot \left[(-\Delta)^{-1/2} \theta\right]$ is an isomorphism from $\dot H^{1}$ into $\dot H^{1}\times \dot H^{1}$ and that  $\theta \longrightarrow \tilde \theta=(1-\alpha^2 \Delta) \theta$ is a bounded operator from $\dot H^{1}$ into $\dot H^{-1}$, and  the Poincar\'e inequality (\ref{PI}), we have the following norm equivalences:
\[||\nabla {\bf v}||_{L^2} \sim~ ||\nabla \theta||_{L^2}\sim~ ||\theta||_{\dot H^1}\sim~||\tilde \theta||_{\dot H^{-1}}.\]

Therefore,
\[ ||F(\tilde \theta_1) -F(\tilde \theta_2)||_{\dot H^{-1} }\le C \left( ||\tilde \theta_1-\tilde\theta_2||_{\dot H^{-1}} \right)\left(||\tilde \theta_1||_{\dot H^{-1}}+||\tilde\theta_2||_{\dot H^{-1}} \right).\]
Consequently, the functional differential equation (\ref{fde}) has short time existence and uniqueness about $t=0$.

Suppose $[0,T_*)$ is the maximal positive interval of existence such that $\tilde \theta \in \mathcal{C}^1([0,T_*),\dot H^{-1}(\Omega))$.

To show the global existence for (\ref{fde}), it is enough to show that the norm $||\tilde \theta||_{\dot H^{-1}}$ stays bounded on the maximal interval of existence. Indeed,
we have  on $[0,T_*)$\[\langle \frac{d}{dt}\tilde \theta,\theta\rangle = \frac{1}{2} \frac{d}{dt} \int(|\theta|^2 + \alpha^2 |\nabla \theta| ^2) \,d{\bf x},\]
and by virtue of Corollary \ref{co1},  Equation (\ref{eq5}) implies that
%
\[ \frac{d}{dt} \int_{\Omega}\left(\theta^2({\bf x},t) + \alpha^2|\nabla \theta({\bf x},t)|^2\right) \,d{\bf x}=0.\]
That is the energy $E(t)$ defined in (\ref{eqe}) is indeed conserved. Therefore
\begin{equation}
\int_{\Omega} \left(\theta^2 + \alpha^2 |\nabla \theta|^2 \right)\,d{\bf x} =\int_{\Omega} \left(\theta_0^2 + \alpha^2 |\nabla \theta_0|^2 \right)\,d{\bf x} \le ||\theta_0||^2_{H^1}\label{Eest}, \text{ for all } \alpha\in (0,1].
\end{equation}
This entrains that the $L^2$ norms of both $\theta$ and its gradient remain bounded. This means that the $\dot H^1$ norm of $\theta$ is bounded or equivalently the norm of $\tilde \theta$ in $\dot H^{-1}$ is bounded.

 A similar argument holds for the negative time interval. This concludes the proof of Theorem \ref{theorem1}.

{\em Remark 1:}\\
It follows from the energy estimate in (\ref{Eest}) above that we have the following bounds:
\begin{equation}\label{E2.1} ||\nabla \theta||_{L^2} \le \frac{1}{\alpha} ||\theta_0||_{\dot H^1} \end{equation}
and \begin{equation} ||\theta||_{L^2} \le  ||\theta_0||_{\dot H^1}.\label{E2}\end{equation}
Therefore in case the solution for the original problem in (\ref{sqg}) develops a singularity in finite time, and if this singular weak solution  is the limit of the   regularized solution, $\theta^{\alpha}$, when $\alpha \longrightarrow 0$, then we at most expect \begin{equation*} ||\nabla \theta^{\alpha}||_{L^2} = O(\frac{1}{\alpha}).\end{equation*}

\section{Higher regularity and a maximum principle}
In this section we discuss  the higher regularity and prove a maximum principle for the regularized problem in (\ref{pb}). We start by proving the following regularity result. The idea of the proof is similar to the presentation in \cite{MB} for the Euler equations.
\begin{theorem}\label{theorem2} ({\bf Regularity):}\\
Let $\theta_0\in H^{m}(\Omega)$, $m\ge 1$, then the solution for the regularized problem (\ref{pb}) $\theta(t)\in \mathcal{C}^1[(-\infty,+\infty),{\dot H}^m]$ (or  $\tilde\theta(t)\in \mathcal{C}^1[(-\infty,+\infty),{\dot H}^{m-2}]$.)
\end{theorem}

{\em Proof:}
\noindent
The case $m=1$ follows from Theorem 1. For $m>1$, we proceed as in the proof of  Theorem 1. We first show local existence and uniqueness in   
$H^m$ then we prove that $||\theta||_{H^m}$ norm remains finite for any  finite interval of time.

It is easy to see that if $m\ge 2$ then \[\theta \in \dot H^m \iff {\bf v } \in \dot H^m\times\dot H^m\iff \tilde \theta \in \dot H^{m-2}(\Omega)\] which implies that \[ \dv{{\bf v}\theta} \in {\dot H}^{m-2}. \] Indeed,
by applying the Gagliardo-Nirenberg and Ladyzhenskaya interpolation inequality, as in Lemma 1, we have
\begin{align}
\int_{\Omega} | D^{m-2} \dv{{\bf v}\theta}|^2 \,d{\bf x} &= \int_{\Omega} | D^{m-2}( {\bf v} \cdot \nabla \theta)|^2 \,d{\bf x} \notag \\ &\le C \sum_{k=0}^{m-2} \int_{\Omega} |D^k{\bf v}\cdot   \nabla D^{m-2-k}\theta|^2 \,d{\bf x} \notag \\
 & \le C \sum_{k=0}^{m-2} ||D^k{\bf v}||^2_{L^4(\Omega)}||\nabla D^{m-k-2}\theta||^2_{L^4(\Omega)}\notag\\
&\le C  \sum_{k=0}^{m-2} ||D^k{\bf v}||_{L^2(\Omega)} ||\nabla D^k{\bf v}||_{L^2(\Omega)}||\nabla D^{m-k-2}\theta||_{L^2(\Omega)}||\nabla \nabla D^{m-k-2}\theta||_{L^2(\Omega)}\notag\\
&\le ||{\bf
v}||_{H^m(\Omega)}^2||{\theta}||_{H^m(\Omega)}^2.\notag\end{align}
Moreover, a similar procedure as in the proof of Theorem 1 applied
to the functional differential equation \[
\frac{d\tilde\theta}{dt}=F(\tilde\theta)\equiv -\dv{{\bf v}\theta}\]
 leads to
\begin{align} ||F(\tilde\theta_1) - F(\tilde  \theta_2)||_{\dot H^{m-2}} &\le  ||\dv{({\bf v}_1-{\bf v}_2)\theta_1}||_{\dot H^{m-2}} + ||\dv{{\bf v}_2(\theta_1-\theta_2)}||_{\dot H^{m-2}}\notag \\
&\le C||\theta_1||_{\dot H^m} ||{\bf v}_1-{\bf v}_2||_{\dot H^m} + ||\theta_1-\theta_2||_{\dot H^m} ||{\bf v}_2||_{\dot H^m}. \notag\end{align}
This completes the proof of short time existence and uniqueness for the equation (\ref{fde}) in $\dot H^{m-2}$.

To show global existence in $ H^m$, it suffices to prove that \[ \Phi(t) = ||D^{m-1}\theta(t)||_{L^2(\Omega)}^2 + \alpha^2 || \nabla D^{m-1}\theta(t)||_{L^2(\Omega)}^2 \]
remains bounded in any finite interval of time. We proceed without proof using the mathematical induction. The case $m=1$ is provided by Theorem \ref{theorem1}.
Assume by induction that $\theta \in C^1[(-\infty,+\infty), H^{m-1}\cap \dot H^1]$.

If $\theta\in \dot H^m\cap\dot H^1$, then $D^{m-1} \theta\in  \dot H^1$. Thus $ \Delta D^{m-1}\theta \in \dot H^{-1}$ and we can write
\begin{align} &\la\frac{\partial}{\partial t}\left(D^{m-1}\theta -\alpha^2 \Delta D^{m-1}\theta\right), D^{m-1}\theta\ra = - \la D^{m-1}({\bf v}\cdot \nabla \theta),D^{m-1}\theta \ra\notag. \\
&\frac{d}{dt} \Phi(t) =-\int_{\Omega} \sum_{k=0}^{m-1} C^{m-1}_k D^k{\bf v} \cdot \nabla (D^{m-k-1}\theta) D^{m-1}\theta \,d{\bf x} \notag \\
&= -\int_{\Omega} \sum_{k=1}^{m-2} C^{m-1}_k D^k{\bf v} \cdot \nabla (D^{m-k-1}\theta) D^{m-1}\theta \,d{\bf x} \notag + \int_{\Omega}   {\bf v} \cdot \nabla (D^{m-1}\theta) D^{m-1}\theta \,d{\bf x}\notag \\ \notag
&=  -\int_{\Omega} \sum_{k=1}^{m-2} C^{m-1}_k D^k{\bf v} \cdot \nabla (D^{m-k-1}\theta) D^{m-1}\theta \,d{\bf x} \notag + \frac12\int_{\Omega}   \dv{{\bf v}  (D^{m-1}\theta)^2   }\,d{\bf x}\notag \\
&=  -\int_{\Omega} \sum_{k=1}^{m-2} C^{m-1}_k D^k{\bf v} \cdot \nabla (D^{m-k-1}\theta) D^{m-1}\theta \,d{\bf x} \notag \\
&\le C\sum_{k=1}^{m-2}||D^k{\bf v} ||_{L^2} || D^{m-1-k}\nabla \theta||_{L^4} ||D^{m-1}\theta||_{L^4}\notag\\
&\le  C\sum_{k=1}^{m-2}||{\bf v} ||_{H^k} || D^{m-1-k}\nabla \theta||^{1/2}_{L^2} || \nabla D^{m-1-k}\nabla \theta||^{1/2}_{L^2} ||D^{m-1}\theta||^{1/2}_{L^2}||\nabla D^{m-1}\theta||^{1/2}_{L^2}\notag\\
&\le C \left(\sum_{k=1}^{m-2}||{\bf v}||_{H^k})\right)||\theta||_{H^{m-1}}||\nabla\theta||_{H^{m-1}} \notag \\
&\le C(\alpha)\left(\sum_{k=1}^{m-2}||{\bf v}||_{H^k}\right) \Phi(t):=\Psi(t) \Phi(t).\notag
\end{align}
Here $\Psi(t) = C(\alpha)\displaystyle\left(\sum_{k=1}^{m-2}||{\bf v}||_{H^k}\right) $. 
Therefore, by using Gronwall's lemma we obtain, \[ \Phi(t) \le \Phi(0) \exp\left(\int_0^t \Psi(s) ds\right),\] which remains bounded on any finite interval of time, since $\Psi(t)$ is bounded by the induction assumption. This completes the proof of the regularity theorem.

Now, we prove that the regularized problem obeys a maximum principle as it is expected that any ``good" regularization of the SQG equation should preserve the physical properties of the original equation. In fact, it is well known   (we can easily show, e.g. let $\alpha=0$ in the proof below) that smooth solutions of the SQG equations obey a maximum principle.

\begin{theorem}
Let the initial condition $\theta_0\in\dot H^1(\Omega)\cap L^{\infty}(\Omega)$. Then the solution $\theta$ of the regularized problem (\ref{pb}) satisfies
\begin{equation}
||\theta(t)||_{L^{\infty}} \le|| \theta_0||_{L^{\infty}} \label{mp}.
\end{equation}
Moreover, if $\theta_0({\bf x}) \ge 0, \forall {\bf x} \in \Omega$ then
\[ \theta(x,t) \ge 0, t>0 ,   {\bf x} \in \Omega.\]\label{thm3}
\end{theorem}

{\em Proof:}\\
Recall the  regularized  solution $\theta$ satisfies the evolution equation
\begin{align}
&\dif{ \theta}{t} - \alpha^2 \Delta \dif{\theta}t+  {\bf v}\cdot \nabla{\theta }= 0 & \text{ in } \Omega \times [0,T]\notag\\
&(-\Delta)^{1/2} \psi = \notag \theta  \\
&\nabla^\bot \psi = {\bf v} ; \int_{\Omega}\theta \,d{\bf x}=0 ; \int_{\Omega}\psi\,d{\bf x}=0 \notag\\
& \theta(x,0) = \theta_0(x) & \text{ in } \Omega \notag .
\end{align}
Notice that $\left(\theta-||\theta_0||_{L^{\infty}}\right)$ satisfies
\[
\dif{ }{t} (\theta-||\theta_0||_{L^{\infty}})- \alpha^2 \Delta \dif{}t(\theta -||\theta_0||_{L^{\infty}})+  {{\bf v}\cdot \nabla(\theta-||\theta_0||_{L^{\infty}}) }= 0.\]
Now we multiple the above evolution equation by $(\theta -||\theta_0||_{L^{\infty}})_+$. Here $w_+ = \max\{w,0\}.$

Observe that if $w\in \dot H^1(\Omega)$ then $w_+\in \dot H^1(\Omega)$ and \[\nabla w_+ = \left\{\begin{array}{cc}\nabla w &\text{ if } w> 0 \\
 0 &\text{ if } w\le 0.\end{array} \right.\]
Therefore \[ \frac{1}{2}\frac{d}{dt}\left[\int_{\Omega} (\theta- ||\theta_0||_{L^{\infty}})_+^2\, d{\bf x} + \alpha^2\int_{\Omega} |\nabla(\theta- ||\theta_0||_{L^{\infty}})_+|^2\, d{\bf x} \right]=0.\]
This yields
\[||(\theta - ||\theta_0||_{L^{\infty}})_+||_{L^2} +\alpha^2 ||\nabla(\theta- ||\theta_0||_{L^{\infty}})_+||^2_{L^2} = ||(\theta(0)-\theta_0||_{L^{\infty}})_+||_{L^2} +\alpha^2 ||\nabla(\theta(0)- ||\theta_0||_{L^{\infty}})_+||^2_{L^2} \]
But the right hand side is zero, because \[ (\theta_0 -||\theta_0||_{L^{\infty}})_+\equiv 0.\]
Thus \[ \left(\theta(t) -||\theta_0||_{L^{\infty}} \right)_+=0, \forall t>0,\]
which implies  \[ \theta({\bf x},t)\le|| \theta_0||_{L^{\infty}}.\]
By using a similar idea, namely by considering the  evolution equation for $(\theta + ||\theta_0||_{L^{\infty}})_{-}$ where $w_-=(-w)_+$ we can show that
 \[\theta({\bf x},t) \ge -||\theta_0||_{L^{\infty}}.\]
Hence (\ref{mp}) follows.

It remains to show the last statement of Theorem \ref{thm3}. Assume $\theta_0({\bf x})\ge 0,{\bf x}\in \Omega$. Multiplying the evolution equation for $\theta$ by $\theta_-=\max\{-\theta,0\}$ and integrating over the domain in a similar fashion as above, yields
\[ \frac{1}{2}\frac{d}{dt}\left[\int_{\Omega} \theta_-^2\, d{\bf x} + \alpha^2\int_{\Omega} |\nabla(\theta_-)|^2\, d{\bf x} \right]=0.\]
i.e.\[ ||\theta_-(t)||^2_{L^2}+\alpha^2||\nabla\theta_-(t)||^2_{L^2} =||\theta_-(0)||^2_{L^2}+\alpha^2||\nabla\theta_-(0)||^2_{L^2} =0\]
because $\theta_0\ge 0$. Thus $\theta_-(t)\equiv 0,\; \forall t>0$. 

\section{Weak convergence and  conditions for blow up for the SQG}
In this section we prove the convergence of the solution of the regularized problem (\ref{pb}) to a weak solution of the original SQG equation (\ref{sqg}), as $\alpha$ tends to zero. We also show that if the original problem has a regular (smooth) solution then the solution for the regularized problem converges strongly to this regular solution. 
Moreover, we show that when these two results are combined together with the energy estimate obtained in section 2.1,  a necessary and sufficient condition for the solution of the original SQG equation to blow up on a finite time interval $[0,T)$ is that the gradient of the regularized solution satisfies
\begin{equation}\label{blw}\sup_{[0,T)}\liminf_{\alpha\goesto 0^+}\alpha^2||\nabla \theta^{\alpha}\nabla ||^2_{L^2} = \epsilon>0 .\end{equation}

Next, the weak convergence of the solution of the regularized problem to a weak solution of the original problem (\ref{sqg}) is discussed.

\begin{theorem}\label{theorem3}Let $T>0$ fixed. Then, the  set of solutions $\theta^{\alpha}, 0<\alpha\le1$  for the regularized problem (\ref{pb}) with   initial condition $\theta_0\in \dot H^1$ is weakly compact in $L^2\big(\Omega\times(-T,T)\big)$. Moreover, if a subsequence of $\theta^{\alpha},\alpha>0$ converges weakly  in $L^2\big(\Omega\times(-T,T)\big)$ to $\bar\theta\in L^2\big(\Omega\times(-T,T)\big)$ when $\alpha\longrightarrow 0$, then $\bar \theta$ is a weak solution for the SQG equations (\ref{sqg}).
\end{theorem}

\noindent
The weak compactness  follows directly from the energy estimate in (\ref{Eest}).
It remains to prove that if $\theta_{\alpha} \longrightarrow \bar \theta$ weakly in $L^2(\Omega)$ then $\bar \theta$ is a weak solution for the SQG equation in (\ref{sqg}).  For this purpose we use the following lemma due to Constantin et al. \cite{canswu}.

\begin{lemma}\label{lemma2}
Let $T>0$ be fixed.
The nonlinear map $\theta \longrightarrow B(\theta,\theta) = \left(\nabla^{\bot}(-\Delta)^{-1/2}\theta \right)\cdot \nabla \theta$  is  weakly continuous  on $L^{2}\big(\Omega\big)$.
\end{lemma}
{\em Proof:} \\
See Appendix B of Constantin et al. \cite{canswu} and references therein.

\noindent
{\em Proof of Theorem \ref{theorem3}:} \\
Let $\phi(x,t)$ be a smooth test function.  Let $\theta^{\alpha}$ be a sequence of solutions for the regularized  SQG equations (\ref{pb}) weakly convergent in $L^{2}(\Omega\times(-T,T))$ to some limit $\bar\theta$.  We have

\[ - \int_{-T}^{T} \int_{\Omega}\theta^{\alpha} \phi_t \,d{\bf x} dt  - \alpha^2 \int_{-T}^{T} \int_{\Omega} \nabla \theta^{\alpha} \nabla \phi_t \,d{\bf x} dt = - \int_{-T}^{T}   \int_{\Omega} B(\theta^{\alpha},\theta^{\alpha}) \phi \,d{\bf x} dt\]
Observe that \[
\left |\alpha^2 \int_{-T}^{T} \int_{\Omega} \nabla \theta \nabla \phi_t \,d{\bf x} dt \right | \le \int_{-T}^{T}\alpha^2 ||\nabla \theta (t)||_{L^2(\Omega)} ||\nabla \phi_t (t)||_{L^2(\Omega)} dt\le C_{\phi} \alpha ||\theta_0||_{\dot H^1}\longrightarrow 0, \text{ when } \alpha \goesto 0.\]
Here we used the energy conservation property derived in the proof of Theorem 1 and the  first upper bound estimate (\ref{E2.1}) in Remark 1. It remains to show that
\begin{equation}\label{r2s}\int_{-T}^{T}   \int_{\Omega} B(\theta^{\alpha},\theta^{\alpha}) \phi \,d{\bf x} dt \goesto \int_{-T}^{T}   \int_{\Omega} B(\bar\theta,\bar\theta) \phi \,d{\bf x} dt.\end{equation}
  Without loss of generality we can assume that the test function $\phi$ is on the form
\[\phi = \psi(t)e^{i{\bf k}\cdot {\bf x}}.\]
It is shown in  the Appendix B of \cite{canswu} that the nonlinearity $B$ satisfies
\begin{equation} \left |\left| (-\Delta)^{-1}\left(B(\theta_1,\theta_2) - B(\theta_2,\theta_2)\right)\right|\right|_w \le C ||\theta_1-\theta_2||_w \left(1+\log\left(1 + ||\theta_1-\theta_2||_w^{-1}\right)\right)\left(||\theta_1||_{L^{2}(\Omega)} + ||\theta_2||_{L^2(\Omega)}\right),\label{appb}\end{equation}
where $||.||_w$ is the weak norm given by \[  ||\theta||_w = \sup_{{\bf j}\in \mathbb{Z}^2\backslash\{0\}}|\hat \theta({\bf j})|.\]
Let \[(-\Delta)^{-1}B(\theta^{\alpha}(t),\theta^{\alpha}(t))= \sum_{{\bf k}\in \mathbb{Z}^2\backslash\{0\}}\hat b_{\alpha,k} (t)e^{i {\bf k}\cdot {\bf x}}\text{ and }(-\Delta)^{-1}B(\bar\theta(t),\bar\theta(t))= \sum_{{\bf k}\in \mathbb{Z}^2\backslash\{0\}}\hat b_k (t)e^{i {\bf k}\cdot {\bf x}}.\]
We have \begin{align}\int_{-T}^{T}   \int_{\Omega} \left(B(\theta^{\alpha},\theta^{\alpha}) -B(\bar\theta ,\bar\theta)\right)\phi \,d{\bf x} dt&=\int_{-T}^{T}   \int_{\Omega} (-\Delta)^{-1}\left(B(\theta^{\alpha} ,\theta^{\alpha}) -B(\bar\theta ,\bar\theta)\right) \Delta \phi \,d{\bf x} dt \notag\\
&= |{\bf k}|^2 \int_{-T}^T\left (\hat b_{\alpha,k} - \hat b_k\right)\psi(t) dt. \end{align}
Combining (\ref{appb}), (\ref{E2}), and the fact that $||\theta^{\alpha}||_w\le ||\theta^{\alpha}||_{L^2}$, we have \[\left|\hat b_{\alpha,k}(t) - \hat b_k(t)\right| \le C(||\theta_0||_{\dot H^1},||\bar\theta (t)||_{L^2}).\]
Therefore, (\ref{r2s}) follows from Lemma \ref{lemma2} and the dominated  convergence theorem of Lebesgue. Thus
\[ - \int_{-T}^{T} \int_{\Omega}\bar\theta\phi_t \,d{\bf x} dt  = - \int_{-T}^{T}   \int_{\Omega} B(\bar\theta,\bar\theta) \phi \,d{\bf x} dt,\] i.e.  the weak limit $\bar\theta$ is a weak solution for the original SQG equation in (\ref{sqg}) on $[-T,T]\times \Omega$.

\noindent
{\em Remark 2}\\
We make the following important observation which provides a sufficient condition for the limiting weak solution for the original SQG equation to blow up in finite time.

Recall the energy conservation property.
\[ \int_{\Omega}\left[ (\theta^{\alpha}(t))^2 +\alpha^2 |\nabla \theta^{\alpha (t)}|^2\right] \,d{\bf x} =  \int_{\Omega}\left[ \theta_0^2+ \alpha^2 |\nabla \theta_0|^2\right] \,d{\bf x}. \]
If \[ \sup_{[0,T)}{\liminf_{\alpha\goesto 0}} \int_{\Omega} \alpha^2 |\nabla \theta^{\alpha}(t)|^2 \,d{\bf x} = \epsilon >0 \]
then either
\begin{itemize}
\item[(i)] $\theta^{\alpha}$ does not converge in norm (i.e. does not converge strongly) to $\bar\theta$ in $L^2(\Omega)$, or
 \item[(ii)] $\bar\theta$ does not conserve energy, i.e.
 \[ \int_{\Omega} \bar\theta^2 \,d{\bf x} \ne \int_{\Omega} \theta_0^2\,d{\bf x}.\]
\end{itemize}
Notice, however, that the weak limit $\bar \theta$ obeys the stability condition
\[ ||\bar\theta||_{L^2}\le ||\theta_0||_{L^2}. \]
Indeed, \[ 0\le ||\theta^{\alpha}-\bar\theta||_{L^2}^2 = ||\theta^{\alpha}||^2_{L^2}-2\la\theta^{\alpha},\bar\theta\ra + ||\bar \theta||^2_{L^2}\le ||\theta_0||^2_{L^2}+\alpha^2||\nabla\theta_0||^2_{L^2} -2\la\theta^{\alpha},\bar\theta\ra + ||\bar \theta||^2_{L^2}.\] When $\alpha\goesto 0$ this yields (because of the weak convergence) \[0\le ||\theta_0||^2_{L^2}-||\bar \theta||^2_{L^2}.\]


Now we prove the following strong convergence theorem, for the regularized problem to the strong solution of the original SQG equation when this latter exists and is regular enough. This in turn  guarantees that we have blowup of the SQG solution if and only if the weak limit in Theorem \ref{theorem3} blows up in finite time.

\begin{theorem}\label{theorem4}
Let the initial condition $\theta_0$ be in $H^2(\Omega)\cap \dot H^1(\Omega)$. Let $\bar\theta\in H^2\cap\dot H^1 $ be
a regular solution for the original SQG equations in (\ref{sqg}) on a finite time interval $[-T,T], T>0$. Then the solution $ \theta^\alpha$ of the regularized problem (\ref{pb}) converges strongly to $\bar \theta$, when $\alpha\goesto 0$. More precisely, we have
\[ \lim_{\alpha\goesto 0} ||(\bar\theta(t)-\theta^{\alpha}(t))||^2_{L^2} +\alpha^2 ||\nabla(\bar\theta(t)-\theta^{\alpha}(t))||^2_{L^2} =0, \text{ uniformly in } [-T,T].\]
 \end{theorem}

\noindent
{\em Proof:}\\
For simplicity in exposition we restrict the discussion to $[0,T]$.
Let $\bar \theta({\bf x},t)\in C^1[[0,T], H^2\cap\dot H^1]$ be a strong solution of the original SQG equation in (\ref{sqg}) with the given initial data and let $\theta^\alpha\in H^2 \cap\dot H^1$ be the corresponding solution for the regularized problem (\ref{pb}).
We have
\begin{equation}\label{exp}\frac{\partial}{\partial t} (\bar \theta -\theta^\alpha) +\alpha^2\frac{\partial}{\partial t} \Delta\theta^{\alpha}+ \dv{ ({\bf \bar v - v^{\alpha}})\bar \theta} - \dv{({\bf \bar v - v^{\alpha}})(\bar \theta-\theta^{\alpha})}+\dv{{\bf \bar v}(\bar \theta-\theta^{\alpha})} = 0,\end{equation} at least in $L^2$, according to the proof of Theorem \ref{theorem1} and the regularity Theorem {\ref{theorem2}}.

First note that
\[  \int_{\Omega} \left(\dv{({\bf \bar v - v}^{\alpha}(\bar \theta-\theta^{\alpha})}\right)(\bar\theta-\theta^{\alpha})\,d{\bf x} =
\int_{\Omega}\left(\dv{{\bf \bar v}(\bar \theta-\theta^{\alpha})}\right)(\bar\theta-\theta^{\alpha})\,d{\bf x}=0. \]
Therefore, the action of the expression in (\ref{exp}) on  $(\bar\theta-\theta^{\alpha})$ yields
\begin{equation*} \frac{1}2 \frac{d}{dt} \int_{\Omega} (\bar \theta -\theta^{\alpha})^2\,d{\bf x} - \alpha^2\int_{\Omega}\left( \Delta \bar\theta_t -\Delta\theta^{\alpha}_t\right)~(\bar \theta -\theta^{\alpha})\,d{\bf x}    + \alpha^2\int_{\Omega} \Delta \bar\theta_t (\bar \theta -\theta^{\alpha})\,d{\bf x}  + \int_{\Omega} [{({\bf \bar v - v^{\alpha}})
 \cdot\nabla \bar \theta}] (\bar \theta - \theta^{\alpha})\,d{\bf x}=0. \end{equation*}

or
\begin{equation*} \frac{1}2 \frac{d}{dt} \int_{\Omega} (\bar \theta -\theta^{\alpha})^2\,d{\bf x} + \alpha^2 \frac{1}{2} \frac{d}{dt}\int_{\Omega} \left| \nabla(\bar \theta -\theta^{\alpha})\right|^2\,d{\bf x}-\alpha^2\int_{\Omega} \nabla \bar\theta_t\cdot \nabla(\bar \theta -\theta^{\alpha})\,d{\bf x}  + \int_{\Omega} [{({\bf \bar v - v^{\alpha}})
 \cdot\nabla \bar \theta}] (\bar \theta - \theta^{\alpha})\,d{\bf x}=0. \end{equation*}
i.e.
\begin{align} &\frac{1}2 \frac{d}{dt} \int_{\Omega} (\bar \theta -\theta^{\alpha})^2\,d{\bf x} + \alpha^2 \frac{1}{2} \frac d{dt}\int_{\Omega} \left| \nabla(\bar \theta -\theta^{\alpha})\right|^2\,d{\bf x}\notag\\& ~~~~~~~~+ \alpha^2\int_{\Omega} \nabla \dv{\bar{\bf v} \bar\theta}\cdot \nabla(\bar \theta -\theta^{\alpha})\,d{\bf x}  + \int_{\Omega} [({\bf \bar v - v^{\alpha}})
 \cdot \nabla \bar \theta] (\bar \theta - \theta^{\alpha})\,d{\bf x}=0. \label{eq23}\end{align}


which implies
\begin{align} &\frac{1}2 \frac{d}{dt} \int_{\Omega} (\bar \theta -\theta^{\alpha})^2 + \alpha^2 \left| \nabla(\bar \theta -\theta^{\alpha})\right|^2\,d{\bf x}\notag\\&~~~~~~~~~\le  \alpha^2|| \nabla \dv{\bar{\bf v} \bar\theta}||_{L^2} ||\nabla(\bar \theta
-\theta^{\alpha})||_{L^2}  + ||\nabla \bar \theta||_{L^{\infty}}||({\bf \bar v - v}^{\alpha})
 ||_{L^2}|| (\bar \theta - \theta^{\alpha})||_{L^2}. \notag\end{align}
But \[\alpha^2|| \nabla \dv{\bar{\bf v} \bar\theta}||_{L^2}
||\nabla(\bar \theta -\theta^{\alpha})||_{L^2} \le \alpha^2|| \nabla
\dv{\bar{\bf v} \bar\theta}||_{L^2}\left( ||\nabla
\bar\theta||_{L^2} + ||\nabla\theta^{\alpha})||_{L^2}
\right)\le\alpha^2
C_{\bar\theta}\left(1+\frac{1}{\alpha}\right)\equiv
\epsilon(\alpha),
\]
where the upper bound estimate $||\nabla \theta^{\alpha}||_{L^2} \le
\frac{1}{\alpha} ||\theta_0||_{\dot H^1}$ in Remark 1 is used.
Therefore,
\begin{align}&
\frac{1}2 \frac{d}{dt} \int_{\Omega} \left[\big(\bar \theta({\bf x},t) -\theta^{\alpha}({\bf x},t)\big)^2 + \alpha^2 \left| \nabla\big(\bar \theta ({\bf x},t)-\theta^{\alpha}({\bf x},t)\big)\right|^2\right]\,d{\bf x}  \notag\\
&\notag~~~~~~~~\le \epsilon(\alpha) + C||\nabla \bar \theta({\bf x},t)||_{L^{\infty}}|| (\bar \theta({\bf x},t) - \theta^{\alpha}({\bf x},t))||^2_{L^2}\notag\\
&\notag~~~~~~~~\le \epsilon(\alpha) + C||\nabla \bar \theta (t)||_{L^{\infty}}\left( || (\bar \theta (t)- \theta^{\alpha}(t))||^2_{L^2} +\alpha^2  || \nabla(\bar \theta(t) - \theta^{\alpha}(t))||^2_{L^2} \right)
\end{align}
which yields
\begin{align} &\int_{\Omega} (\bar \theta({\bf x},t) -\theta^{\alpha}({\bf x},t))^2 + \alpha^2 \left| \nabla(\bar \theta ({\bf x},t)-\theta^{\alpha}({\bf x},t))\right|^2\,d{\bf x}\notag \\ &~~~~~~~\le \epsilon(\alpha)T+ \int_0^t C||\nabla \bar \theta (s)||_{L^{\infty}}\left( || (\bar \theta(s) - \theta^{\alpha}(s))||^2_{L^2} +\alpha^2  || \nabla(\bar \theta(s) - \theta^{\alpha}(s))||^2_{L^2}\right) ds.\end{align}
Therefore, by using Gr\"onwall's lemma, we have
\[ ||(\bar\theta(t)-\theta^{\alpha}(t))||^2_{L^2} +\alpha^2 ||\nabla(\bar\theta(t)-\theta^{\alpha}(t))||^2_{L^2} \le T \epsilon(\alpha) \exp\left(\int_0^T C||\nabla \bar \theta (s)||_{L^{\infty}}\,ds\right)\goesto 0, \text{ when }\alpha\goesto 0.\]


\noindent
Finally, we  show that the condition anticipated  in (\ref{blw}) is indeed necessary and sufficient for the original problem to have a singular solution. Therefore providing a systematic  and practical procedure, relying only on the behavior of the regularized problem (\ref{pb}),  for  detecting  the eventual blowup in finite time of the smooth solutions for the SQG equations in (\ref{sqg}). More precisely we have the following result.
\begin{theorem}
Let $[0,T^*)$ be the maximal interval of existence for the
regular solution $\bar\theta$ for the original SQG problem. That is, it develops a singularity in its gradient at time $t=T^*$  if and only if
\begin{equation}\sup_{[0,T^*)}\liminf_{\alpha\goesto 0^+}\alpha^2||\nabla \theta^{\alpha}  ||^2_{L^2} = \epsilon>0 \label{nc}.\end{equation}
\end{theorem}
{\em Proof:}
First we show that it is a sufficient condition.
Assume that a sequence of solutions $\theta^{\alpha}$ for the regularized problem (\ref{pb}) converges weakly to a weak solution $\bar \theta\in L^2(\Omega)$ for the original problem (\ref{sqg}).
As stated in Remark 2, if \[\sup_{[0,T^*)}\liminf_{\alpha\goesto 0^+}\alpha^2||\nabla \theta^{\alpha} ||^2_{L^2} = \epsilon>0,\]
then either $\theta^{\alpha}$ does not converge strongly to $\bar \theta$ or that $\bar \theta$ is not a regular solution.
Theorem \ref{theorem4}, however, guarantees that a regular solution is necessarily a strong limit of $\theta^{\alpha}$. Therefore (\ref{nc}) is a sufficient condition for blow up in finite time.

Now assume that $\bar \theta\in \dot H^1(\Omega)$ is a regular solution for the SQG equations on a maximal interval of existence $[0,T^*)$ such that
\[ \limsup_{t\goesto T^*}||\nabla \bar\theta||_{L^2} =+\infty.\]
According to Theorem \ref{theorem4}, for $0\le t<T^*$ fixed, we have \[\lim_{\alpha \goesto 0} \alpha||\nabla \theta^{\alpha}(t)-\nabla\bar\theta(t)||_{L^2}=0.\]
Let $t$ be sufficiently close to $T^*$ so that $||\nabla \bar\theta(t)||_{L^2}>\frac{\delta}{\alpha}$ where $\delta>0$ is fixed and $\alpha>0$ sufficiently small so that $\alpha ||\nabla\theta^{\alpha}(t)-\nabla\bar\theta(t)||_{L^2} <\delta/2$, we have
\[ \alpha||\nabla\theta^{\alpha}(t)||_{L^2} \ge\alpha ||\nabla\bar\theta(t)||_{L^2} -\alpha ||\nabla\theta^{\alpha}(t)-\nabla\bar\theta(t)||_{L^2}> \frac{\delta}{2} .\]
Therefore \[\limsup_{t\goesto T^*}\liminf_{\alpha\goesto0} \alpha||\nabla\theta^{\alpha}||_{L^2}= \epsilon>0.\]

\section*{Acknowledgment}
The work of B.K. is partly supported by a grant from the National Sciences and Engineering Research Council of Canada.
The work of E.S.T. was supported in part by the NSF grant
no.~DMS-0504619, the BSF grant no.~2004271,  and the ISF grant no.~120/06.


 \end{document}